\documentclass[graybox]{svmult}

\usepackage{mathptmx}       
\usepackage{helvet}         
\usepackage{courier}        
\usepackage{makeidx}         
\usepackage{graphicx}        
\usepackage{multicol}        
\usepackage[bottom]{footmisc}

\makeindex             

\begin{document}

\title*{Understanding the Dynamics of Collision and Near-Collision Motions in the $N$-Body Problem}
\titlerunning{Collisions and Near-Collisions} 
\author{Lennard F. Bakker}
\institute{Lennard F. Bakker \at Department of Mathematics, Brigham Young University, USA \email{bakker@math.byu.edu}}

\maketitle

\abstract
{Although rare, collisions of two or more bodies in the $N$-body problem are apparent obstacles at which Newton's Law of Gravity ceases to make sense. Without understanding the nature of collisions, a complete understanding of the N-body problem can not be achieved. Historically, several methods have been developed to probe the nature of collisions in the $N$-body problem. One of these methods removes, or regularizes, certain types of collisions entirely, thereby relating not only analytically, but also numerically, the dynamics of such collision motions with their near collision motions. By understanding the dynamics of collision motions in the regularized setting, a better understanding of the dynamics of near-collision motions is achieved.}

\section{Introduction}
\label{sec:1}

For ages, humankind has observed the regular and predicable motion of the planets and other bodies in the solar system and asked, will the motion of the bodies in the solar system continue for ever as they are currently observed? This philosophical question is the object of the mathematical notion of stability. A difficulty in applying the notion of stability to the motion of the solar system is that of collision and near-collision motions of bodies in the solar system. Collision and near-collision motions do occur in the solar system. Section \ref{sec:2} recounts a few of these that have been observed or predicted.

The standard mathematical model for understanding the motion of planets and other bodies in the solar system is the Newtonian $N$-Body problem, presented in Section \ref{sec:3}. In included here are some of the basic features and mathematical theory of the Newtonian $N$-Body Problem, its integrals or constants of motion, special solutions such as periodic solutions, and the notions of stability and linear stability of periodic solutions and their relationship.

The notions and basic theory of collisions and singularities in the Newtonian $N$-Body Problem is presented in Section \ref{sec:4}. This includes a discussion of the probabilities of collisions, and the regularization or the lack thereof for collisions. A collision motion is rare in that is has a probability of zero of occurring, whereas a near collision motion has a positive probability of occurring. Regularization is a mathematical technique that removes the collision singularities from the Newtonian $N$-Body Problem and enables an analysis of near-collision motions in terms of collision motions through the continuous dependence of motions on initial conditions. This regularization is illustrated in the collinear $2$-Body Problem, the simplest of all the $N$-Body Problems.

Recent results are presented in Section \ref{sec:5} on the analytic and numerical existence and numerical stability and linear stability of periodic motions with regularizable collisions in various $N$-Body Problems with $N=3$ and $N=4$. Although fictitious, these periodic motions with regularizable collisions provide a view of their near-collision motions which could be motions of the bodies in the $N$-Body Problem that are collision-free and bounded for all time.

\section{Phenomenon}
\label{sec:2}

Collisions and near-collisions of two or more solar system bodies are apparent obstacles at which Newton's Law of Gravity becomes problematic. Velocities of colliding bodies become infinite at the moment of collision, while velocities of near-colliding bodies become very large as they pass by each other. Both of these situations present problems for numerical estimates of the motion of such bodies.

Although collisions are rare, historical evidence of collisions of solar system bodies is viewable on the surface of the Earth and the Moon \cite{Ce}. Only recently have collisions of solar systems bodies actually been observed. As the comet Shoemaker-Levy 9 approached Jupiter it was torn apart into fragments by tidal forces. In July of 1994, at least 21 discernible fragments of Shoemaker-Levy 9 collided with Jupiter. These were the {\it first} ever observed collisions of solar system bodies. An animation of some of the fragments of Shoemaker-Levy 9 colliding with Jupiter can be found at www2.jpl.nasa.gov/sl9/anim.html. 

Near-collision motion are less rare than collisions. As of March 2012, there are nearly 9000 known near-Earth asteroids\footnote{See http://neo.jpl.nasa.gov/stats/}, of which 1306 are potentially hazardous to Earth\footnote{See http://neo.jpl.nasa.gov/neo/groups.html}. One of these potentially hazardous asteroids, named 2012 DA14, was discovered in 2012. This asteroid will pass by Earth on February 15, 2013, coming closer to the Earth than satellites in geostationary orbit\footnote{See article about 2012 DA14 posted March 6, 2012 on MSNBC.com}. How close will 2012 DA14 pass by Earth? A mere 17000 miles (27000 km)\footnote{See article about 2012 DA14 posted March 8, 2012 on Earthsky.org}. In cosmic terms, this close shave of 2012 DA14 with Earth in 2013 is a near-collision motion.

\section{The $N$-Body Problem}
\label{sec:3}

To model collision and near-collision motions we make some simplifying assumptions and use Newton's Inverse Square Law of Gravity. We assume that all the bodies are idealized as particles with zero volume (i.e., as points), that no particle is torn apart by tidal forces, that the mass of each particle never changes, and that besides Newton's Law of Gravity there are no other forces acting on the bodies. Under these assumptions we would think of Shoemaker-Levy 9 as not being torn apart by tidal forces, but as colliding with Jupiter as a whole. 

\subsection{The Equations}
\label{sec:3.1}

The particles modeling the bodies move in three-dimensional Euclidean space which we denote by ${\mathbf R}^3$. For a positive integer $N\geq 2$, suppose there are $N$ particles with positions ${\mathbf q}_j\in{\mathbf R}^3$ and masses $m_j>0$, $j=1,\dots,N$. The distance between two of the particles is denoted by
\[ r_{jk} = \vert {\mathbf q}_j-{\mathbf q}_k\vert, \ \ j\ne k,\]
which is the standard Euclidean distance between two points in ${\mathbf R}^3$. The Newtonian $N$-Body Problem is the system of second-order nonlinear differential equations
\[ m_j {\mathbf q}_j^{\prime\prime} = \sum_{k\ne j} \frac{Gm_jm_k({\mathbf q}_k-{\mathbf q}_j)}{r_{jk}^3},\ \ j=1,\dots,N,\]
where ${}^\prime = d/dt$ for a time variable $t$ and $G=6.6732\times 10^{-11}\ {\rm m}^2/{\rm s}^2{\rm kg}$. By an appropriate choice of units of the ${\mathbf q}_j$, we will assume that $G=1$ because we are investigating the qualitative or geometric, rather than the quantitative, behavior of collision and near-collision motions.

By the standard existence, uniqueness, and extension theory in differential equations (see \cite{Ch} for example), the initial value problem
\begin{equation}\label{IVP} m_j {\mathbf q}_j^{\prime\prime}  = \sum_{k\ne j} \frac{m_jm_k({\mathbf q}_k-{\mathbf q}_j)}{r_{jk}^3}, \ \  {\mathbf q}_j(t_0)  = {\mathbf q}_j^0, \ \
{\mathbf q}_j^{\prime}(t_0)  = {\mathbf  q}_j^{\prime 0}, \end{equation}
has a unique solution
\[ {\mathbf q}(t) = ( {\mathbf q}_1(t),\dots,{\mathbf q}_N(t))\]
defined on a maximal interval of definition $(t^-,t^+)$ as long as $r_{jk}\ne 0$ for all $j\ne k$ at $t=t_0$. Such a solution ${\mathbf q}(t)$ describes a motion of the $N$ particles.

Not every initial value problem (\ref{IVP}) will have a solution ${\mathbf q}(t)$ with $t^-=-\infty$ and $t^+=\infty$. A solution with either $t^->-\infty$ or $t^+<\infty$ experiences a singularity at the finite endpoint of its maximal interval of definition. The notion of a singularity is addressed in Section \ref{sec:4.1}.

\subsection{Integrals}
\label{sec:3.2}

An integral of motion of the Newtonian $N$-Body Problem is a differentiable function $F$ of the position ${\mathbf q}$ and/or the velocity ${\mathbf q}^\prime$ and/or the masses ${\mathbf m}=(m_1,\dots,m_N)$ such that
\[ \frac{d}{dt} F({\mathbf q}(t),{\mathbf q}^\prime(t),{\mathbf m}) = 0, \ \ t\in (t^-,t^+).\]
Along a solution ${\mathbf q}(t)$ an integral $F$ of motion satisfies
\[ F({\mathbf q}(t),{\mathbf q}^\prime(t),{\mathbf m}) = F({\mathbf q}(t_0),{\mathbf q}^\prime(t_0),{\mathbf m}),\ \ t^-<t<t^+,\]
i.e., the value of $F$ is constant along the solution. The Newtonian $N$-Body Problem has ten known integrals of motion.

The translation invariance of the equations of the Newtonian $N$-Body Problem give rise to $6$ integrals of motion. With $M=\sum_{j=1}^N m_j$, three of these are given by the components of the center of mass vector
\[ {\mathbf C} = \frac{1}{M} \sum_{j=1}^m m_j {\mathbf q}_j,\]
and three more are given by the components of the linear momentum vector
\[ {\mathbf L} = \frac{1}{M} \sum_{j=1}^N m_j{\mathbf q}_j^\prime.\]
Typically, both of these are set to ${\mathbf 0}$ so that the {\it relative motion} of the $N$ particles is emphasized.

The rotational symmetry of the equations of the Newtonian $N$-Body Problem give rise to $3$ more integrals of motion. These integrals are given by the components of the angular momentum vector,
\[ {\mathbf A} = \sum_{j=1}^N m_j {\mathbf q}_j\times {\mathbf q}^\prime_j.\]
The angular momentum plays a key role in understanding collisions in the $N$-Body Problem, as we will see later on.

There is one more integral of motion of the Newtonian $N$-Body Problem. The {\it self-potential} (or negative of the potential energy) is
\[ U = \sum_{j<k} \frac{m_j m_k}{r_{ik}^2}.\]
The {\it kinetic energy} is
\[ K = \frac{1}{2} \sum_{j=1}^N m_j {\mathbf q}^\prime_j\cdot{\mathbf q}^\prime_j.\]
The {\it total energy}
\[ H = K-U\]
is an integral of motion for the Newtonian $N$-Body Problem.

In the late 1800's, the mathematical strategy for ``solving'' the Newtonian $N$-Body Problem was to find enough ``independent'' integrals of motion \cite{Sa}. This would implicitly give each solution as the curve of intersection of the hypersurfaces corresponding to the integrals of motion. Each solution ${\mathbf q}(t)$ is a curve in ${\mathbf R}^{6N}$. However, the intersection of the hypersurfaces of the ten integrals of motion gives a $6N-10>1$ dimension hypersurface in ${\mathbf R}^{6N}$, which is not a curve! The ten known integrals of motion are independent of each other (one is not a function of the others), and are algebraic functions of positions, velocities, and masses. Are there any more algebraic integrals of motion? This was answered a long time ago in 1887-1888 by Bruns \cite{Br}.

\begin{theorem} There are no algebraic integrals of motion independent of the ten known integrals of motion.
\end{theorem}

\noindent Consequently, new integrals of motion, if any, can not be algebraic!  In 1893, Newcomb \cite{Ne} lamented that no additional integrals had been found to enable the implicit solution of the $3$-Body Problem. It is well-known that the Newtonian $2$-Body Problem can be solved implicitly\footnote{See en.wikipedia.org/wiki/Gravitational\_two-body\_problem}, but all attempts to solve the $N$-Body Problem with $N\geq 3$ have been futile\footnote{Karl Sundman did solve the $3$-Body Problem when ${\mathbf A}\ne 0$ by convergent power series defined for all time, but the series converge too slowly to be of any theoretic or numerical use \cite{Sa}.}.

Typically then the solution ${\mathbf q}(t)$ of the initial value problem (\ref{IVP}) is estimated numerically. From the constant total energy $H$ along a solution ${\mathbf q}(t)$, we observe that if any of the distances $r_{jk}$ get close to $0$, i.e., at least two of the particles are near collision, the self-potential becomes large and the kinetic energy becomes large too. The latter implies that the velocity of at least one of the particles becomes large, and the linear momentum ${\mathbf L}$ along ${\mathbf q}(t)$ implies that the velocity of at least two particles becomes large. In particular, from the equations of the Newtonian $N$-Body Problem, the particles that are near collision are the one with the large velocities. These large velocities presents problems for the numerical estimates of such a solution.

\subsection{Special Solutions}
\label{sec:3.3}

Rather than solving the $N$-Body Problem for all of its solutions by finding enough independent integrals of motion, it is better to examine special solutions with particular features. The simplest solutions to find are equilibrium solutions, where the position ${\mathbf q}_j(t)$ of each particle is constant for all time. But the Newtonian $N$-Body Problem has none of these (see p.29 in \cite{MHO}). The next simplest solutions are periodic solutions, i.e., there exist $T>0$ such that ${\mathbf q}(t+T) = {\mathbf q}(t)$ for all $t\in{\mathbf R}$. These are part of the larger collection of solutions ${\mathbf q}(t)$ with $t^-=-\infty$ and $t^+=\infty$ that are bounded. Such solutions must have a particular total energy (see p. 160 in \cite{Sa}).

\begin{theorem} If a solution ${\mathbf q}(t)$ of the Newtonian $N$-Body Problem exists for all time and is bounded, then the total energy $H<0$.
\end{theorem}

\noindent Consequently, any periodic solution ${\mathbf q}(t)$ of the Newtonian $N$-Body Problem must has negative total energy. This is why in the search for periodic solutions the total energy is always assigned a negative value.

\subsection{Stability}
\label{sec:3.4}

A periodic solution ${\mathbf q}(t)$ of the Newtonian $N$-Body Problem gives a predictable future: we know with certainty what the positions of the $N$ particles will be at any time $t>0$. But what if our measurements of the initial conditions ${\mathbf q}(0)$ and ${\mathbf q}^\prime(0)$ are slightly off? A solution $\tilde{\mathbf q}(t)$ with initial conditions near ${\mathbf q}(0)$ and ${\mathbf q}^\prime(0)$ will stay close to ${\mathbf q}(t)$ for a time, by a property of solutions of initial value problems called continuity of solutions with respect to initial conditions (see \cite{Ch}). But if it stays close for all $t>0$, we think of ${\mathbf q}(t)$ as being stable.

To {\it quantify} this notion of stability for a periodic solution, we use a Poincar\'e section which is a hyperplane $S$ containing the point $({\mathbf q}(0),{\mathbf q}^\prime(0))$ that is transverse to the curve $({\mathbf q}(t),{\mathbf q}^\prime(t))$. If ${\mathbf x} = (\tilde{\mathbf q}(0),\tilde{\mathbf q}^\prime(0))$ is a point on $S$ near the $({\mathbf q}(0),{\mathbf q}^\prime(0))$, then $P({\mathbf x})$ is next point where the the curve $(\tilde{\mathbf q}(t),\tilde{\mathbf q}^\prime(t))$ intersects $S$\footnote{For an illustration of this see en.wikipedia.org/wiki/Poincar\'e\_map}, and $P^2({\mathbf x})$ is the next point, and so on. The initial condition ${\mathbf x}^0 = ({\mathbf q}(0),{\mathbf q}^\prime(0))$ is a {\it fixed point} of this Poincar\'e map $P$ from $S$ to $S$, i.e., $P({\mathbf x}^0)={\mathbf x}^0$.

\begin{definition}\label{stability} The periodic solution ${\mathbf q}(t)$ is stable if for every real $\epsilon>0$ there exist a real $\delta>0$ such that $\vert P^k({\mathbf x}) - {\mathbf x}^0\vert<\epsilon$ for all $k=1,2,3,\dots,$ whenever $\vert {\mathbf x} - {\mathbf x}^0\vert<\delta$.
\end{definition}

\noindent When ${\mathbf q}(t)$ is not stable, there are solutions which start nearby but eventually move away from ${\mathbf q}(t)$, and we say that ${\mathbf q}(t)$ is {\it unstable}.

Showing directly that ${\mathbf q}(t)$ is stable or unstable is very difficult. Instead, the related concept of linearized stability is investigated, at least numerically. The derivative of the Poincar\'e map at the fixed point ${\mathbf x}^0$ is a square matrix $DP({\mathbf x}^0)$.

\begin{definition}\label{linearstability} A periodic solution ${\mathbf q}(t)$ is 
\begin{enumerate}
\item spectrally stable\footnote{There is a more restrictive notion of spectral stability known as linear stability that requires additional technical conditions on the square matrix $DP({\mathbf x}^0)$.} if all the eigenvalues of $DP({\mathbf x}^0)$ have modulus one, and is
\item linearly unstable if any eigenvalue of $DP({\mathbf x}^0)$ has modulus bigger than one.
\end{enumerate}
\end{definition}

\noindent In 1907, Liapunov \cite{Li} established a connection between the stability of Definition \ref{stability} and the linearized stability of Definition \ref{linearstability}.

\begin{theorem} If a periodic solution ${\mathbf q}(t)$ is stable then it is spectrally stable, and if ${\mathbf q}(t)$ is linearly unstable, then it is unstable.
\end{theorem}

\noindent If a periodic solution is shown numerically to be linearly unstable, then by Theorem 3 the periodic solution is unstable. On the other hand, if a periodic solution is shown numerically to be spectrally stable, it may be stable or unstable. Examples exist with spectrally stable fixed points of maps like $P$ that are unstable (see \cite{SM}).

The notion of stability for a non-periodic solution, such as the motion of the sun and planets in the solar system, is harder to grasp. Here is a sampling of the history and opinions on this stability problem. In 1891, Poincar\'e commented that the stability of the solar system had at that time already preoccupied much time and attention of researchers (see p. 147 in \cite{DH}). In 1971, Siegel and Moser lamented that a resolution of the stability problem for the $N$-Body Problem would probably be in the distant future (see p. 219 in \cite{SM}). In 1978, Moser noted that the answer to the stability of the solar system was still not known (see p. 127 in \cite{DH}). In 2005, Saari stated that a still unresolved problem for the $N$-Body Problem is that of stability (see p. 132 in \cite{Sa}). Meyer, Hall, and Offin commented how little is known about the stability problem and how difficult it was to get (see p. 229 in \cite{MHO}).

In 1996, Diacu and Holmes suggested that the solar system should be considered stable (in a weak sense) if no collisions occur among the sun and the planets, and no planet ever escape from the solar system (see. p.129 in \cite{DH}). In this weak sense of stability, the solar system is stable for the next few billion years according to numerical work of Hayes \cite{Ha} in 2007. Much longer-term numerical studies of the solar system by Batygin and Laughlin  \cite{BL} in 2008 using small changes in the initial conditions suggest that Mercury could fall into the sun in 1.261Gyr\footnote{Gyr means giga-year or 1,000,000,000 years}, or that Mercury and Venus could collide in 862Myr\footnote{Myr means mega-year or 1,000,000 years} and Mars could escape from the solar system in 822Myr. The Newtonian $N$-Body Problem thus suggests that in the near future, the Solar System should be free of collisions of planets and the Sun, with no planets escaping the solar system. But this still leaves open the possibility that smaller objects, such as asteroids and comets, could collide with any of the planets in the short and long term. Recall that there are nearly 9000 of those near-Earth asteroids to consider, with 2012 DA14 making its near-collision approach with Earth on February 15 of 2013.

\section{Collisions}
\label{sec:4}

Either in the short term of the long term, collisions put a {\it wrench} into the question of any notion of stability. Why should a solution or any nearby solution of the Newtonian $N$-Body Problem be defined for all time? Remember that Shoemaker-Levy 9 has $t^+<\infty$!

\subsection{Singularities}
\label{sec:4.1}

Collisions are one of the {\it two} kinds of singularities in the Newtonian $N$-Body Problem. The solution ${\mathbf q}(t)$ of initial value problem (\ref{IVP}) is real analytic (i.e., a convergent power series) on an interval $(t_0-\delta,t_0+\delta)$ for some $\delta>0$, as long as $r_{jk}\ne 0$ for all $j\ne k$ at $t_0$. By a process called analytic continuation (see for example \cite{MaH}), the interval $(t_0-\delta,t_0+\delta)$ can be extended to the maximal interval $(t^-,t^+)$.

\begin{definition}
A {\it singularity} of the Newtonian $N$-Body Problem is a time $t=t^+$ or $t^-$ when $t^+<\infty$ or $t^->-\infty$.
\end{definition}

In 1897, Painlev\'e \cite{Pa} characterized a singularity of the Newtonian $N$-Body Problem, using the quantity
\[ r_{\rm min}(t) = \min_{j\ne k} r_{jk}(t)\]
determined by a solution ${\mathbf q}(t)$.

\begin{theorem} A singularity for the Newtonian $N$-body Problem occurs at time $t=t^*$ if and only if $r_{\rm min}(t)\to 0$ as $t\to t^*$.
\end{theorem}

\noindent An understanding what this means is obtained by considering the {\it collision set}
\[ \Delta = \bigcup_{j\ne k} \{ {\mathbf q} : {\mathbf q}_j={\mathbf q}_k\}\subset ({\mathbf R}^3)^N,\]
which is the set of points where two or more of the $N$-particles occupy the same position. Painlev\'e's characterization means that ${\mathbf q}(t)$ approaches the collison set, i.e.,
\[ {\mathbf q}(t)\to \Delta\ \ {\rm as}\ \ t\to t^*\]
when $t^*$ is a singularity of the Newtonian $N$-Body Problem. Painlev\'e's chararterization introduces two classes of singularities.

\begin{definition} A singularity $t^*$ of the Newtonian $N$-Body Problem is a collision singularity when $q(t)$ approaches a specific point of $\Delta$ as $t\to t^*$. Otherwise the singularity $t^*$ is a non-collision singularity.
\end{definition}

Only collision singularities can occur in the Newtonian $2$-Body Problem because it can be implicitly solved. In 1897, Painlev\'e \cite{Pa} showed that only one other Newtonian $N$-Body problem has only collision singularities.

\begin{theorem} In the $3$-Body Problem, all singularities are collision singularities. 
\end{theorem}

\noindent Unable to extend his result to more than $3$ bodies, Painlev\'e conjectured that there exist non-collision singularities in the Newtonian $4$ or larger Body Problem. In 1992, Xia \cite{Xi} mostly confirmed Painlev\'e's conjecture, giving an example in the Newtonian $5$-Body Problem.

\begin{theorem} There exist non-collision singularities in the $N$-Body Problem for $N\geq 5$.
\end{theorem}

\noindent That leaves unresolved the question of the existence of non-collision singularities in the Newtonian $4$-Body Problem.

An understanding of what a non-collision singular looks like is obtained by considering one-half of the {\it polar moment of inertia} of the Newtonian $N$-Body Problem:
\[ I = \frac{1}{2} \sum_{j=1}^N m_j {\mathbf q}_j\cdot{\mathbf q}_j.\]
This scalar quantity measures the ``diameter'' of the $N$ particles in the Newtonian $N$-Body Problem. In 1908, von Zeipel \cite{Ze} characterized a collision singularity in terms of the polar moment of inertia.

\begin{theorem}\label{vonZeipel} A singularity of the Newtonian $N$-Body Problem at $t=t^*$ is a collision if and only if $I$ is bounded as $t\to t^*$.
\end{theorem}

\noindent This implies that for a non-collision singularity, at least one of the $N$-particles has to achieve an infinite distance from the origin in just a finite time. This is a rather strange thing for Newton's Law of Gravity to predict. On the other hand, by Theorem \ref{vonZeipel}, for a collision singularity, all of the positions of the $N$ particles remain bounded at the moment of the singularity.

A total collapse is an example of a collision singularity in the $N$-Body Problem for which all $N$ particles collide at the same point at the singularity $t^*$. For a solution ${\mathbf q}(t)$, the quantity \[ r_{\rm max} = \max_{j\ne k} r_{jk}(t)\]
characterizes a total collapse: a total collapse occurs at $t^*$ if and only if
\[ r_{\rm max}(t) \to 0 {\rm\ as\ } t\to t^*.\]
There is a relationship between total collapse and the angular momentum that was known by Weierstrass and established by Sundman (see \cite{Sa}).

\begin{theorem}\label{angular} If ${\mathbf A}\ne 0$, then $r_{\rm max}(t)$ is bounded away from zero.
\end{theorem}

\noindent This does not preclude the collision of less than $N$ particles when ${\mathbf A}\ne 0$, as will be illustrated for certain Newtonian $N$-Body Problems in Section \ref{sec:5}.

\subsection{Improbability}
\label{sec:4.2}

Recall that there are 1306 potentially hazardous near-Earth asteroids. What are the chances that Earth will be hit by a near-Earth asteroid, or Jupiter will be hit by another comet? Well, it depends on the arrangement of the particles.

\begin{definition} A solution ${\mathbf q}(t)$ is called collinear if the $N$ particles always move on the same fixed line in ${\mathbf R}^3$. Otherwise it is called non-collinear.
\end{definition}

\noindent Every collinear solution has zero angular momentum because ${\mathbf q}_j(t)$ is parallel with ${\mathbf q}_j^\prime(t)$ for all $t\in (t^-,t^+)$. In 1971 and 1973, Saari \cite{Sa71,Sa73} established the probability of collisions.

\begin{theorem}\label{collisionprobability}The probability that a non-collinear solution ${\mathbf q}(t)$ will have a collision is zero. Every collinear solution ${\mathbf q}(t)$ has a collision.
\end{theorem}

With collision singularities being rare for a non-collinear $N$-Body Problem, why bother to study them? Diacu and Holmes (see p. 84 and p. 103 in \cite{DH}) argue for the study of collision singularities because without such a study, a complete understanding of the Newtonian $N$-Body Problem could not be achieved. In particular, solutions near collision singularities could behave strangely, and the probability of a solution coming close to a collision singularity is positive and thus can not be ignored. Understanding then the collision singularities enables an understanding of the near-collision solutions.

\subsection{Regularization}
\label{sec:4.3}

Regularization is one method by which we can get an understanding of a collision singularity. To {\it regularize} a collision means to extend the solution beyond the collision through an elastic bounce without loss or gain of total energy in such a way that all of the solutions nearby have continuity with respect to initial conditions, i.e., they look like the extended collision solution for a time (see p. 104 and p. 107 in \cite{DH}). Regularization is typically done by a Levi-Civita type change of the dependent variables, and a Sundman type change of the independent variable (see \cite{Ce}),
that together {\it removes} the collision singularity from the equations. We illustrate this regularization in the simplest of the $N$-Body Problems.

In the Collinear $2$-Body Problem (or Col2BP for short), the positions of the two particles are the scalar quantities $q_1$ and $q_2$. If $x=q_2-q_1$ is the distance between the particle with mass $m_1$ at $q_1$ and the particle with mass $m_2$ at $q_2>q_1$, then the Col2BP takes the form
\begin{equation}\label{Col2BP} x^{\prime\prime} = - \frac{m_1+m_2}{x^2}, \ x>0,\end{equation}
and the total energy takes the form
\begin{equation}\label{totalenergy} H = \frac{m_1m_2}{2(m_1+m_2)} (x^\prime)^2 - \frac{m_1m_2}{x}.\end{equation}
As $x\to 0$ the two particles approach collision, and the total energy implies that the two particles collide with an infinite velocity,
\[ (x^\prime)^2\to\infty.\]
To regularize the {\it binary collision} (or total collapse) in this problem, define a new independent variable $s$ and a new dependent variable $w$ by
\[ \frac{ds}{dt} = \frac{1}{x}, \ \ w^2 = x,\]
where the former is the Sundman type change of the independent variable, and the latter is the Levi-Civita type change of the dependent variable. If $\ \dot{}= d/ds$, the second-order equation (\ref{Col2BP}) becomes
\begin{equation}\label{regularizedODE} w^2\big[ 2w\ddot w - 2\dot w^2 + (m_1+m_2)\big]=0,\end{equation}
and the total energy (\ref{totalenergy}) becomes
\begin{equation}\label{Renergy} Hw^2 = \frac{2m_1m_2}{m_1+m_2}\dot w^2 - m_1m_2.\end{equation}
As $w\to 0$, the second-order equation (\ref{regularizedODE}) makes sense (no dividing by zero), and the total energy (\ref{Renergy}) implies that
\[ (\dot w)^2 \to \frac{m_1+m_2}{2},\]
which is a finite nonzero velocity! The collision singularity has been regularized.

The regularized nonlinear second-order equation (\ref{regularizedODE}) can actually be solved! Solving the total energy (\ref{Renergy}) for $2(\dot w)^2$ and substituting this into the second-order equation (\ref{regularizedODE}) gives
\begin{equation}\label{linearODE} 2w^3\left[ \ddot w - \frac{(m_1+m_2)H}{2m_1m_2} w \right] = 0.\end{equation}
This makes sense when $w=0$, i.e., the moment of collision! For negative $H$, the linear second-order equation\footnote{This is a simple harmonic oscillator for $H<0$ whose solutions are in terms of cosine and sine.} inside the square brackets in (\ref{linearODE}) solves to give a real analytic stable periodic solution $w(s)$ which experiences a collision every half period in terms of the regularized time variable $s$. The corresponding solution $x(t)$ is periodic and experiences a collision once a period in terms of the original time variable $t$. This doubling of the number of collisions per period is because the change of dependent variable $w^2=x$ has $w(s)$ ``doubling'' $x(t)$ in that $w(s)$ passes through $0$ twice a period, going from positive to negative and then negative to positive, while $x(t)$ is positive except at collision where it is zero.

The binary collision singularity in the Newtonian $2$-Body Problem can be regularized in a similar but more complicated way than what was done above for the Col2BP (see \cite{Sa}). By Theorem \ref{angular}, a solution of the $2$-Body Problem with nonzero angular momentum does not experience a collision or total collapse. A nonzero angular momentum near-collision solution looks like the zero angular momentum collision solution\footnote{Binary Star Systems are known to exist in the Universe. The Newtonian $2$-Body Problem predicts stability for a Binary Star System, a collision-free solution that is bounded for all time.}.  The regularized $2$-Body Problem provides good numerical estimates of the motion because there are no infinite velocities!

\subsection{McGehee}
\label{sec:4.4}

What about regularization of a triple collision, when three of the particles meet? In 1974, McGehee \cite{McG} showed that regularization of a triple collision is in general not possible\footnote{This is achieved by ``blowing-up'' the triple collision singularity and slowing down the motion as the particles approach a triple collision. This setting does allow for good numerical estimates of near-triple collisions.}. Starting close together, two solutions that approach a near triple collision can describe {\it radically different} motions after the near triple collision. This kind of behavior is known as ``sensitive dependence on initial conditions,'' and is an antithesis of stability. Triple collisions present a numerical nightmare! By extension, collisions with four or more particles present the same nightmare! So the only regularizable collisions are those that are essentially a binary collision.

\section{Results}
\label{sec:5}

Spectrally stable periodic solutions have been found in Newtonian $N$-Body Problems with regularizable collisions for $N\geq 3$. Three of these situations discussed here are the Collinear $3$-Body Problem (or Col3BP), the Collinear Symmetric $4$-Body Problem (or ColS4BP), and the Planar Pairwise Symmetric $4$-Body Problem (or PPS4BP). There are other Newtonian $N$-Body Problems where periodic solutions with regularizable collisions whose existence has been given analytically \cite{Ya1,Ya2,Sh}, some of whose stability (in the sense of Definition \ref{stability}) and linear stability (as defined in Definition \ref{linearstability}) has been numerically determined \cite{BS,Wa,Ya1,Ya2}.

\subsection{Col3BP}
\label{sec:5.1}

As a subproblem of the Newtonian $3$-Body Problem, the Col3BP requires that the three particles always lie on the same line through the origin. The positions of the three particles in the Col3BP are the scalars $q_1$, $q_2$, and $q_3$ which can be assumed to satisfy
\[ q_1\leq q_2\leq q_3.\]
By Theorem \ref{collisionprobability}, collisions always occur in the Col3BP. Because the three particles are collinear for all time, their angular is zero, and by Theorem \ref{angular} a total collapse is possible\footnote{Initial conditions leading to total collapse in the equal mass Col3BP are easy to realize: set $q_1=-1$, $q_2=0$, and $q_3=1$ with the initial velocity of each particle set to $0$.} in the Col3BP. In 1974, S.J. Aareth and Zare \cite{AZ} showed that any two of the three possible binary collisions in the $3$-Body Problem are regularizable\footnote{A good numerical model for the Sun-Jupiter-Shoemaker-Levy 9 or Earth-Moon-2012DA14 situation is regularized $3$-Body Problem of Aarseth and Zare.}. In 1993, Hietarinta and Mikkola \cite{HM} used Aarseth and Zare's regularization \cite{AZ} to regularize the binary collisions $q_1=q_2$ and $q_2=q_3$ in the Col3BP.

In 1956, Schubart \cite{Sc} numerically found a periodic orbit in the equal mass Col3BP of negative total energy in which the inner particle oscillates between binary collisions with the outer particles. In 1977, H\'enon \cite{He} numerically extended Schubart's periodic solution to arbitrary masses and investigated their linear stability. In 1993, Hietarinta and Mikkola \cite{HM} also numerically investigated the linear stability of Schubart's periodic solution for arbitrary masses. Together they showed that Schubart's periodic solution is spectrally stable for certain masses, and linearly unstable for the remaining masses. Hietarinta and Mikkola \cite{HM} further numerically investigated the Poincar\'e section for Schubart's periodic solution for arbitrary masses, showing when there is stability as described in Definition \ref{stability}.  In 2008, Moeckel \cite{Mo} and Venturelli \cite{Ve} separately proved the analytic existence of Schubart's solution when $m_1=m_3$ and $m_2$ is arbitrary. Only recently, in 2011, did Shibayama \cite{Sh} analytically prove the existence of Schubart's periodic solution for arbitrary masses in the Col3BP.

Schubart's period solution for the Col3BP is also a periodic solution of the $3$-Body Problem, where in the latter the continuity with respect to initial conditions can be see for near-collision solutions. For example, Schubart's periodic solution for the nearly equal masses
\[ m_1=0.333333,\ m_2=0.333334,\ m_3=0.333333\]
is spectrally stable. Considered in $3$-Body Problem, Schubart's periodic solution for these mass values remains spectrally stable \cite{He}, and numerically the near-collision solutions in the Newtonian $3$-Body Problem behave like Schubart's periodic solution. It is therefore possible that in the $3$-Body Problem, there are solutions near Schubart's periodic solution that are free of collisions and bounded for all time. Imagine, as did H\'enon \cite{He}, of Newton's Law of Gravity predicting a triple star system that is free of collisions and bounded for all time!

\subsection{ColS4BP}
\label{sec.5.2}

As a subproblem of the Newtonian $4$-Body Problem, the ColS4BP requires that the four particles always lie on the same line through the origin. The positions of the four particles are the scalars $q_1$, $q_2$, $q_3$, and $q_4$ that satisfy
\[ q_4=-q_1,\ q_3=-q_2,\ q_1\geq 0,\ q_2\geq 0,\]
and
\[ -q_1 \leq -q_2\leq 0\leq q_2\leq q_1\]
with masses
\[ m_1=1,\ m_2=m>0,\ m_3=m,\ m_4=1.\]
The angular momentum for all solutions of the ColS4BP is zero because of the collinearity, and so by Theorem \ref{angular} a total collapse is possible. There are two kinds of non-total collapse collisions in the ColS4BP: the binary collision of the inner pair of particles of mass $m$ each, i..e, $q_2=0$, and the {\it simultaneous} binary collision of the two outer pairs of particles, i.e., $q_1=q_2>0$. In 2002 and 2006, Sweatman \cite{SW1,SW2} showed, by adapting the regularization of Aarseth and Zare \cite{AZ}, that these non-total collapse collisions in the ColS4BP are regularizable.

Sweatman \cite{SW1,SW2} numerically found a Schubart-like periodic solution in the ColS4BP with negative total energy for arbitrary $m$ where the outer pairs collide in a simultaneous binary collision at one moment and then the inner pair collides at another moment. He determined numerically that this Schubart-like periodic solution is spectrally stable when
\[ 0<m<2.83\ {\rm and}\  m>35.4,\]
and is otherwise linearly unstable. In 2010, Bakker et al \cite{BOYSR} verified Sweatman's linear stability for the Schubart-like periodic solution using a different technique. In 2011-2012, Ouyang and Yan \cite{OY2}, Shibayama \cite{Sh}, and Huang \cite{Hu} proved separately the analytic existence of the Schubart-like periodic solution in the ColS4BP.

\subsection{PPS4BP}
\label{sec:5.3}

The PPS4BP has two particles of mass $1$ located at the planar locations
\[ {\mathbf q}_1 {\rm\ and\ }{\mathbf q}_3=-{\mathbf q}_1,\]
and two particles of mass $0<m\leq 1$ located at the planar locations
\[ {\mathbf q}_2 {\rm\ and\ } {\mathbf q}_4=-{\mathbf q}_2.\]
The four particles in the PPS4BP need not be collinear, so that the angular momentum need not be $0$. Unlike the ColS4BP, total collapse can be avoided in the PPS4BP by Theorem \ref{angular} when the angular momentum is not zero. Like the ColS4BP, there are two kinds of non-total collapse collisions in the PPS4BP: simultaneous binary collisions when ${\mathbf q}_1={\mathbf q}_2$ and ${\mathbf q}_3={\mathbf q}_4$, or when ${\mathbf q}_1 = {\mathbf q}_4$ and ${\mathbf q}_2={\mathbf q}_3$; and binary collisions when ${\mathbf q}_1=0$ or when ${\mathbf q}_2=0$. In 2010, Sivasankaran, Steves, and Sweatman \cite{SSS} showed that these non-total collapse collisions in the PPS4BP are regularizable.

The Schubart-like periodic solution in the ColS4BP is also a periodic solution of the PPS4BP, where in the latter the continuity with respect to initial conditions can be observed for near-collision solutions. However, as shown by Sweatman \cite{SW2},  in the PPS4BP the Schubart-like periodic solution of the ColS4BP becomes linearly unstable for
\[ 0<m<0.406, {\rm \ and\ } 0.569<m<1.02\]
as well as $2.83<m<35.4$, while it remains spectrally stable for
\[ 0.407<m<0.567{\rm\ and\ } m>35.4.\]
By long-term numerical integrations for the Schubart-like periodic solution as a solution of the PPS4BP, Sweatman \cite{SW2} showed that stability in the sense of Definition \ref{stability} is possible when $0.407<m<0.567$ and when $m>35.4$. It is therefore possible for these values of $m$ that near Schubart's periodic solution there are collision-free solutions of the PPS4BP that are bounded for all time.

In 2011, adapting the regularization of Aarseth and Zare \cite{AZ} to simultaneous binary collisions, Bakker, Ouyang, Yan, and Simmons \cite{BOYS} proved the analytic existence of a non-collinear periodic solution in the equal mass PPS4BP. This periodic solution has zero angular momentum, negative total energy,  and alternates between a simultaneous binary collision of the symmetric pairs in the first and third quadrant where ${\mathbf q}_1={\mathbf q}_2$ and ${\mathbf q}_3={\mathbf q}_4$, and the simultaneous binary collision of the symmetric pairs in the second and fourth quadrants where ${\mathbf q}_1={\mathbf q}_4$ and ${\mathbf q}_2={\mathbf q}_3$. Bakker, Ouyang, Yan and Simmons \cite{BOYS} then numerically extended this non-collinear periodic simultaneous binary collision solution to unequal masses $0<m<1$. In 2012, Bakker, Mancuso, and Simmons \cite{BMS} have numerically determined that the non-collinear periodic  simultaneous binary collision solution is spectrally stable when
\[ 0.199<m<0.264{\rm\ and\ } 0.538<m\leq 1\]
and is linearly unstable for the remaining values of $m$. Long-term numerical integrations of the regularized equations done by Bakker, Ouyang, Yan, and Simmons \cite{BOYS} suggest instability when $0.199<m<0.264$ and stability when $0.538<m\leq 1$ in the sense of Definition \ref{stability}. For these latter values of $m$ could the near-collision solutions in the PPS4BP that look like the non-collinear periodic simultaneous binary collision solution be collision-free and bounded for all time?

\section{Future Work?}
\label{sec:6}

Both the ColS4BP and the PPS4BP are subproblems of the Newtonian $4$-Body Problem, where the non-total collapse collisions in the former two problems are regularizable. What is not known is how to, if possible, regularize binary collisions and simultaneous binary collisions in the Newtonian $4$-Body Problem within one coordinate system \footnote{During the Special Session on Celestial Mechanics at the American Mathematical Society's Sectional Conference in April 2011 at the College of the Holy Cross, Worcester, Massachusetts, Rick Moeckel put forth the problem of finding an elegant coordinate system for the Newtonian $4$-Body Problem in which regularizes binary collisions and simultaneous binary collisions and blows up all triple collisions and total collapse. The regularization of binary collisions and simultaneous binary collisions can be achieved within multiple coordinate systems, with one coordinate system for each regularizable collision.}. If such a regularization is possible, then all of the periodic solutions thus known in the ColS4BP and PPS4BP would also be periodic solutions of the Newtonian $4$-Body Problem and the investigation of their stability and linear stability in the Newtonian $4$-Body Problem could begin. With more possible perturbations of initial conditions in the Newtonian $4$-Body Problem as compared with the PPS4BP, a loss of spectral stability could indeed happen as it did with going from the ColS4BP to the PPS4BP. But some of the spectral stability might survive passage from the PPS4BP to the Newtonian $4$-Body Problem, giving the possibility of near-collision solutions that are collision-free and bounded for all time.

\begin{acknowledgement} The author expresses appreciation for the referee's comments and feedback that improved the quality of this paper. The author also expresses thanks to the organizers of the year-long seminar series held at Virginia State University.
\end{acknowledgement}


\begin{thebibliography}{99.}
 

\bibitem{AZ} Aarseth, S.J., Zare, K.: A regularization of the three-body problem, Cel. Mech. {\bf 10} (1974)

\bibitem{BOYSR} Bakker, L.F.,  Ouyang, T.,  Yan, D., Simmons, S.C., Roberts, G.E.:  Linear Stability for Some Symmetric Periodic Simultaneous Binary Collision Orbits in the Four-Body Problem, Celest. Mech. Dynam. Astron. {\bf 108} (2010)

\bibitem{BOYS}  Bakker, L.F., Ouyang, T.,  Yan, D., Simmons, S.C.: Existence and Stability of Symmetric Periodic Simultaneous Binary Collision Orbits in the Planar Pairwise Symmetric Four-Body Problem, Celest. Mech. Dynam. Astron. {\bf 110} (2011)

\bibitem{BMS} Bakker, L.F., Mancuso, S.C., Simmons, S.C.: Linear stability analysis of symmetric periodic simultaneously binary collision orbits in the planar pairwise symmetric four-body problem, J. Math. Anal. Appl. {\bf 392} (2012)

\bibitem{BS} Bakker, L.F., Simmons, S.C.: Stability of the Rhomboidal Symmetric-Mass Orbit, submitted to J. Math. Anal. Appl (Aug. 2012), \url{http://arxiv.org/pdf/1208.3183.pdf}

\bibitem{BL} Batygin, K., Laughlin, G.: On the dynamical stability of the Solar System, The Astrophysical Journal {\bf 683} (2008)

\bibitem{Br} Bruns, H.: \"Uber die Integrale des Vielk\"orper-Problems, Acta Math. {\bf 11} (1887-1888)

\bibitem{Ce} Celletti, A.: Singularities, Collisions and Regularization Theory, Lecture Notes in Physics {\bf 590}, Benest D. and Froeschl\'e, C. (eds), Singularities in Gravitational Systems, Springer, New York (2002)

\bibitem{Ch} Chicone, C.: Ordinary Differential Equations with Applications, Texts in Applied Mathematics {\bf 34}, Springer, New York (1999)

\bibitem{DH} Diacu. F., Holmes, P.: Celestial encounters: the origin of chaos and stability, Princeton University Press, Princeton (1996)

\bibitem{Ha} Hayes, W.: Is the Outer System chaotic?, \url{http://arxiv.org/abs/astro-ph/0702179v1}

\bibitem{He} H\'enon, M.: Stability of interplay orbits, Cel. Mech. {\bf 15} (1977).

\bibitem{HM} Hietarinta, J., Mikkola, S.: Chaos in the one-dimensional gravitational three-body problem, Chaos {\bf 3} (1993)

\bibitem{Hu} Huang, Hsin-Yuan: Schubart-like orbits in the Newtonian collinear four-body problem: a variational proof, Dis. Con. Dyn. Sys. {\bf 32} (2012)

\bibitem{Li} Liapunov, A.: Prob\'eme g\'en\'eral de la stabilit\'e du mouvement. Ann. Fac. Sci. Toulouse {\bf 9} 1907

\bibitem{MaH} Marsden, G.E., Hoffman, M.J.: Basic Complex Analysis, Second Edition, W.H. Freeman and Company, New York (1987)

\bibitem{McG} McGehee, R.: Triple collision in the collinear three-body problem, Inventiones Mathematicae {\bf 27} (1974)

\bibitem{MHO} Meyer, K.R., Hall, D.R., Offin, D.: Introduction to Hamiltonian Dynamical Systems and the $N$-Body Problem, Second Edition, Springer, New York (2009)

\bibitem{Mo} Moeckel, R.: A Topological Existence Proof for the Schubart Orbits in the Collinear Three-Body Problem, Dis. Con. Dyn. Syst. Series B {\bf 10} (2008)

\bibitem{Ne} Newcomb, S.: Modern Mathematical Thought. In: Bull of New York Math Soc. {\bf 4} (1893)

\bibitem{OY2} Ouyang, T., Yan, D.: Periodic solutions with alternating singularities in the collinear four-body problem, Celest. Mech. Dynam. Astron.  {\bf 109} (2011)

\bibitem{Pa} Painlev\'e, P.: Lecons sur la th\'eorie analytic de equations diff\'erentielles, Herman, Paris (1897)

\bibitem{Sa71} Saari, D.G.: Improbability of collisions in Newtonian gravitational systems, Trans. Amer. Math. Soc. {\bf 162} (1971)

\bibitem{Sa73} Saari, D.G.: Improbability of collisions in Newtonian gravitational systems II, Trans. Amer. Math. Soc. {\bf 181} (1973)

\bibitem{Sa} Saari, D.G.: Collisions, Rings, and Other Newtonian $N$-Body Problems, CBMS {\bf 104}, American Mathematical Society, Providence, Rhode Island (2005)

\bibitem{Sc} Schubart, j.: Numerische Aufsuchung periodischer L\"osungen im Dreik\"orperproblem, Astronomische Nachriften {\bf 283} (1956)

\bibitem{Sh} Shibayama, M.: Minimizing Periodic Orbits with Regularizable Collisions in the $n$-Body Problem, Arch. Rational Mech. Anal. {\bf 199} (2011)

\bibitem{SM} Siegel, C.L., Moser, J.K.: Lectures on Celestial Mechanics. Springer, New York (1971)

\bibitem{SSS} Sivasankaran, A., Steves, B.A., Sweatman, W.L.: A global regularisation for integrating the Caledonian symmetric four-body problem, Celest. Mech. Dynam. Astron. {\bf 107} (2010)

\bibitem{SW1} Sweatman, W.L.: Symmetrical one-dimensional four-body problem, Celest. Mech. Dynam. Astron. {\bf 82} (2002)

\bibitem{SW2} Sweatman, W.L.: A Family of Symmetrical Schubart-Like Interplay Orbits and their Stability in the One-Dimensional Four-Body Problem, Celest. Mech. Dynam. Astron. {\bf 94} (2006)

\bibitem{Ve} Venturelli, A.: A Variational Proof of the Existence of Von Schubart's Orbit, Dis. Con. Dyn. Syst. Series B, {\bf 10} (2008)

\bibitem{Wa} Waldvogel, J.: The rhomboidal symmetric four-body problem, Celest. Mech. Dyn. Astron. {\bf 113} (2012)

\bibitem{Xi} Xia, Z.: The existence of Noncollision Singularities in Newtonian Systems, Ann. Math. {\bf 135} (1992)

\bibitem{Ya1} Yan, D.: Existence and linear stability of the rhomboidal periodic orbit in the planar equal mass four-body problem, J. Math. Anal. Appl. {\bf 388} (2012).

\bibitem{Ya2} Yan, D.: Existence of the Broucke orbit and its linear stability, J. Math. Anal. Appl. {\bf 389} (2012)

\bibitem{Ze} von Zeipel, E.H.: Sur les singularit\'es du probl\`eme des n corps, Ark. Mat. Astron. Pys. {\bf 4} (1908)


%
%
%
%
%

%
%
%
%
%
%
%
%
%
%
%
%
%
%
%
%
%
%
%
%
%
%
%
%
%
%
%
%
%
%
%
%
%


\end{thebibliography}
\end{document}